\documentclass{article}

\newtheorem{fed}{\textbf{Definition}}[section]

\newtheorem{lemma}[fed]{\textbf{Lemma}}

\newtheorem{rem}[fed]{\textbf{Remark}}
\newtheorem{prop}[fed]{\textbf{Proposition}}
\newtheorem{cor}[fed]{\textbf{Corollary}}
\usepackage{amssymb,bbm,graphicx,epsfig,psfrag,epic,eepic,latexsym}
\usepackage{amsmath}
\usepackage{mathrsfs}
\usepackage[dvips]{color}

\newcommand{\N}{\mathbb{N}}

\newcommand{\R}{\mathbb{R}}

\begin{document}
\title{First steps in the geography of scale Hilbert structures}
\author{Urs Frauenfelder \footnote{Department of Mathematics and
Research Institute of Mathematics, Seoul National University}}
\maketitle

\begin{abstract}
Scale structures were introduced by H.\,Hofer, K.\,Wysocki, and
E.\,Zehnder. In this note we define an invariant for scale Hilbert
spaces modulo scale isomorphism and use it to distinguish large
classes of scale Hilbert spaces.
\end{abstract}

\tableofcontents

\section[Statement of the main results]{Statement of the main results}

Scale structures on a Banach space were introduced by H.\,Hofer,
K.\,Wysocki, and E.\,Zehnder, see \cite{hofer1,
hofer-wysocki-zehnder1, hofer-wysocki-zehnder2}. They observed that
on a scale Banach space a new notion of smoothness can be defined
which still meets the chain rule. Therefore scale structure give
rise to new smooth structures in infinite dimensions. Manifolds
modelled on this new smooth structures are called scale manifolds.
These provide the first step in the construction of polyfolds which
in turn can be used to deal with transversality issues in Symplectic
Field theory, Gromov-Witten theory, or Floer theory. The author's
interest in this new smooth structures in infinite dimensions is
based on the following guiding principle. He believes that the
various Floer homologies should be interpretable as \emph{Morse
homology on scale manifolds}. Such a unified framework would lead to
various simplifications of the existing theory, since gluing and
transversality issues could be referred to the general set-up
currently developed by H.\,Hofer, K.\,Wysocki, and E.\,Zehnder, and
have not be checked anymore individually for each Floer homology.

In this note we introduce a first invariant to distinguish different
Hilbert scale structures and we construct various examples of
nonisomorphic Hilbert scale structures. The restriction to the
Hilbert case instead of the more general Banach case is justified by
our intension to apply scale structures to Floer homology. In Floer
homology one need to have metrics since one has to be able to define
a gradient.

We first recall the definition of a Hilbert scale structure which is
due to H.\,Hofer, K.\,Wysocki, and E.\,Zehnder.
\begin{fed}\label{sca}
A \emph{scale Hilbert space} is a tuple
$$\mathcal{H}=\big\{\big(H_k, \langle \cdot, \cdot \rangle_k\big)\big\}_{k
\in \N_0}$$ where for each $k \in \N_0$ the pair $\big(H_k,\langle
\cdot, \cdot \rangle_k\big)$ is a real Hilbert space and the vector
spaces $H_k$ build a nested sequence
 $H=H_0 \supset H_1 \supset H_2 \supset \ldots$
such that the following two axioms hold.
\begin{description}
 \item[(i)] For each $k \in \N$ the inclusion
 $\big(H_k, \langle \cdot, \cdot \rangle_k\big) \hookrightarrow
 \big(H_{k-1}, \langle \cdot, \cdot \rangle_{k-1}\big)$ is compact.
 \item[(ii)] For each $k \in \N_0$ the subspace
 $H_\infty=\bigcap_{n=0}^\infty H_n$ is dense in $H_k$ with respect
 to the topology induced from $\langle \cdot, \cdot \rangle_k$.
\end{description}
\end{fed}
\begin{rem}\label{scar}
\emph{If $H_0$ is finite dimensional, then the second axiom in
Definition~\ref{sca} implies that $H_k=H_0$ for every $k \in \N_0$.
On the other if $H$ is infinite dimensional, then the first axiom
implies that $H_k \neq H_0$ for every $k \in \N$.}
\end{rem}
We next introduce the notion of isomorphism between two scale
Hilbert spaces. Hence suppose that $\mathcal{H}=\{(H_k,\langle
\cdot, \cdot \rangle_k)\}$ and $\mathcal{H}'=\{(H_k',\langle \cdot,
\cdot \rangle'_k)\}$ are two scale Hilbert spaces. We denote for $k
\in \N_0$ by $||\cdot||_k$ and $||\cdot||'_k$ the norms on $H_k$,
respectively $H'_k$, induced from the scalar products $\langle
\cdot, \cdot \rangle_k$ and $\langle \cdot, \cdot \rangle'_k$.
\begin{fed}
A \emph{scale isomorphism} $\Phi$ from $\mathcal{H}$ to
$\mathcal{H}'$ is a bijective linear map
$$\Phi \colon H_0 \to H'_0$$
satisfying the following two axioms.
\begin{description}
 \item[(i)] For each $k \in \N_0$ the map $\Phi$ restricts to a
 bijection
 $$\Phi_k \colon H_k \to H_k', \quad \Phi_k=\Phi|_{H_k}.$$
 \item[(ii)] For each $k \in \N_0$ there exists a constant $c_k>0$,
 such that
 $$\frac{1}{c_k}||h||_k \leq ||\Phi(h)||'_k \leq c_k ||h||_k,
 \quad \forall\,\,h \in H_k.$$
\end{description}
\end{fed}
\begin{fed}
Two scale Hilbert spaces $\mathcal{H}$ and $\mathcal{H}'$ are called
\emph{scale isomorphic}, if there exists a scale isomorphism from
$\mathcal{H}$ to $\mathcal{H}'$.
\end{fed}
Since on a finite dimensional vector space all scalar products are
equivalent it follows from Remark~\ref{scar} that in each finite
dimension there is precisely one Hilbert scale structure up to scale
isomorphism. We therefore restrict in the following our attention to
Hilbert scale structures in infinite dimensions. We introduce the
following set
$$\mathscr{S}=\big\{\mathcal{H}\,\,\textrm{scale Hilbert space},\,\,
\dim(H_0)=\infty\big\}\big/\sim$$ where the equivalence relation is
given by scale isomorphism. Geography of Hilbert scale structures
refers to the description of the set $\mathscr{S}$.

To construct examples as well as invariants for scale Hilbert spaces
we introduce the notion of a scale Hilbert $n$-tuple for $n \in \N$.
\begin{fed}
A \emph{scale Hilbert $n$-tuple} is a tuple
$$\mathcal{H}=\big\{\big(H_k, \langle \cdot, \cdot
\rangle_k\big)\big\}_{0 \leq k \leq n-1}$$ where for each $k \in
\{0,\ldots, n-1\}$ the pair $\big(H_k,\langle \cdot, \cdot
\rangle_k\big)$ is a real Hilbert space and the vector spaces $H_k$
build a nested sequence
 $H=H_0 \supset H_1 \supset H_2 \supset \ldots \supset H_{n-1}$
such that the following two axioms hold.
\begin{description}
 \item[(i)] For each $k \in \{0,\ldots, n-1\}$ the inclusion
 $\big(H_k, \langle \cdot, \cdot \rangle_k\big) \hookrightarrow
 \big(H_{k-1}, \langle \cdot, \cdot \rangle_{k-1}\big)$ is compact.
 \item[(ii)] For each $k \in \{0,\ldots,n-1\}$ the subspace
 $H_{n-1}$ is dense in $H_k$ with respect
 to the topology induced from $\langle \cdot, \cdot \rangle_k$.
\end{description}
We refer to scale Hilbert $2$-tuples as \emph{scale Hilbert pairs}
and to scale Hilbert $3$-tuples as \emph{scale Hilbert triples}.
\end{fed}
The notion of scale isomorphism between scale Hilbert $n$-tuples is
the same as the one for scale Hilbert spaces and two Hilbert
$n$-tuples are called scale isomorphic if there exists a scale
isomorphism between them. We next introduce for each $n \in \N$ the
set
$$\mathscr{S}_n=\big\{\mathcal{H}\,\,\textrm{scale Hilbert $n$-tuple},\,\,
\dim(H_0)=\infty\big\}\big/\sim.$$ We further denote by
$\widetilde{\mathcal{F}}$ the space of all functions $f \colon \N
\to (0,\infty)$ which are monotone and unbounded. We say that $f_1,
f_2 \in \widetilde{\mathcal{F}}$ are equivalent if there exists
$c>0$ such that
$$\frac{1}{c} f_1(n) \leq f_2(n) \leq c f_1(n), \qquad n \in \N$$
and we write $f_1 \sim f_2$ for equivalent functions. We introduce
the quotient
$$\mathcal{F}=\widetilde{\mathcal{F}}/\sim.$$
By $\ell^2$ we denote as usual the Hilbert space of all square
summable sequences together with its standard inner product. For $f
\in \widetilde{\mathcal{F}}$ we define $\ell^2_f \subset \ell^2$ to
be the vector space of all sequences $x=(x_1,x_2,\ldots)$ satisfying
$$||x||_{\ell^2_f}=\sqrt{\sum_{n=1}^\infty f(n)x_n^2} <\infty.$$
The inner product
$$\langle x,y \rangle_{\ell^2_f}=\sum_{n=1}^\infty f(n)x_n y_n, \quad
x,y \in \ell^2_f$$ endows $\ell^2_f$ with the structure of a Hilbert
space. We are now in position to state our first main result.
\\ \\
\textbf{Theorem~A } \emph{There is a bijection between the sets
$\mathcal{F}$ and $\mathscr{S}_2$ given by the map $[f] \mapsto
[(\ell^2,\ell^2_f)]$.}
\\ \\
As an immediate consequence of Theorem~A we obtain the following
Corollary.
\begin{cor}\label{a1}
Assume that $\mathcal{H}=\{(H_k,\langle \cdot,\cdot \rangle_k)\}$ is
a scale Hilbert space, then for every $k \in \N_0$ the Hilbert space
$(H_k,\langle \cdot, \cdot \rangle_k)$ is separable.
\end{cor}
\begin{rem}
\emph{Since every separable Hilbert space is actually isometric to
$\ell^2$, it follows that in a scale Hilbert space all Hilbert
spaces are isometric to each other and the geography question for
scale Hilbert spaces is reduced to the question how these infinitely
many $\ell^2$-spaces can be nested into each other.}
\end{rem}
Theorem~A can be used to define invariants for scale Hilbert spaces
modulo the equivalence relation given by scale isomorphism. Let
$\Delta \subset \N_0 \times \N_0$ be the set
$$\Delta=\big\{(i,j) \in \N_0 \times \N_0: i<j\big\}.$$
Denote by $\mathfrak{J}_2 \colon \mathcal{F} \to \mathscr{S}_2$ the
bijection $[f] \mapsto [(\ell^2,\ell^2_f)]$ given by Theorem~A. Now
we introduce the map
$$\mathfrak{K} \colon \mathscr{S} \to
\mathrm{Map}(\Delta,\mathcal{F})$$ which is given for an infinite
dimensional scale Hilbert space $\mathcal{H}=\{(H_k,\langle
\cdot,\cdot \rangle_k)\}$ by the formula
\begin{equation}\label{inva}
\mathfrak{K}([\mathcal{H}])(i,j)=
\mathfrak{J}_2^{-1}\Big(\big[\big(H_i,\langle \cdot,\cdot
\rangle_i\big), \big(H_j,\langle \cdot,\cdot
\rangle_j\big)\big]\Big),\quad (i,j) \in \Delta.
\end{equation}
The same kind of invariant can also be used for scale Hilbert
$n$-tuples for every $n \in \N$ satisfying $n \geq 2$. Indeed, let
$\Delta_n \subset \{0,\ldots,n-1\}\times \{0,\ldots,n-1\}$ be the
set
$$\Delta_n=\big\{(i,j) \subset \N_0 \times \N_0: i<j <n\big\}.$$
Then we define the map
$$\mathfrak{K}_n \colon \mathscr{S}_n \to
\mathrm{Map}(\Delta_n,\mathcal{F})$$ by the same formula
(\ref{inva}) as before. That the invariants $\mathfrak{K}$ and
$\mathfrak{K}_n$ are well defined, i.e.~independent of the choice of
the representative $\mathcal{H}$ is a further Corollary of
Theorem~A.
\begin{cor}\label{a3}
The maps $\mathfrak{K}$ and $\mathfrak{K}_n$ for $n=\{2,3,\ldots\}$
are well defined.
\end{cor}
We can furthermore use Theorem~A to construct a class of examples of
scale Hilbert spaces which are not scale isomorphic to each other.
We first introduce the set
$$\widetilde{\mathfrak{F}}=\mathrm{Map}(\N,\widetilde{\mathcal{F}}).$$
We say that $F_1, F_2 \in \widetilde{\mathfrak{F}}$ are equivalent,
if $F_1(k)$ is equivalent to $F_2(k)$ in $\widetilde{\mathcal{F}}$
for every $k \in \N$. We write again $F_1 \sim F_2$ for equivalent
$F_1, F_2 \in \widetilde{\mathfrak{F}}$ and set
$$\mathfrak{F}=\widetilde{\mathfrak{F}}/\sim.$$
For $F \in \widetilde{\mathfrak{F}}$ we introduce the nested
sequence of Hilbert spaces
$$\ell^{2,F}=\ell^{2,F}_0\supset \ell^{2,F}_1 \supset
\ell^{2,F}_2 \supset \ldots$$ where we set
$$\ell^{2,F}_0=\ell^2$$
and for each $k \in \N$
$$\ell^{2,F}_k=\ell^2_{\prod_{j=1}^k F(j)}$$
where the product of two functions $f_1,f_2 \in
\widetilde{\mathcal{F}}$ is defined pointwise as
$$(f_1 \cdot f_2)(\nu)=f_1(\nu)\cdot f_2(\nu), \quad \nu \in \N.$$
An analogous procedure gives us finite nested sequences of Hilbert
spaces. For $n \in \N$ satisfying $n \geq 2$ we put
$$\widetilde{\mathfrak{F}}_n=\mathrm{Map}\Big(\{1,\ldots n-1\},
\widetilde{\mathcal{F}}\Big).$$ Defining the equivalence relation
pointwise as before we set
$$\mathfrak{F}_n=\widetilde{\mathfrak{F}}_n/\sim$$
and for $F \in \widetilde{\mathfrak{F}}_n$ we introduce again the
now finite nested sequence of Hilbert spaces
$$\ell^{2,F}=\ell^{2,F}_0\supset \ell^{2,F}_1 \supset
\ell^{2,F}_2 \supset \ldots \supset \ell^{2,F}_{n-1}.$$ As a second
Corollary of Theorem~A we can draw the following assertion.
\begin{cor}\label{a2}
For each $n \in \N$ satisfying $n \geq 2$ there is an injective map
$$\mathfrak{J}_n \colon \mathfrak{F}_n \to \mathscr{S}_n, \quad
[F] \mapsto [\ell^{2,F}].$$ Moreover, there is an injective map
$\mathfrak{J} \colon \mathfrak{F} \to \mathscr{S}$ given by the same
formula.
\end{cor}
\begin{rem}\label{nisu}
\emph{For $n=2$ the above Corollary is just a special case of
Theorem~A if one uses the canonical identification of
$\mathfrak{F}_2$ with $\mathcal{F}$ given by $[F] \mapsto [F(1)]$.
In this case the map $\mathfrak{J}_2$ is actually surjective. So one
might wonder if this continues to hold for larger $n \in \N$.
However, surjectivity of $\mathfrak{J}_n$ actually already fails for
$n=3$ which is the content of Corollary~\ref{b1}.}
\end{rem}
From Corollary~\ref{a2} we obtain some information about which
values of the invariant $\mathfrak{K}$ are realizable by scale
Hilbert spaces. Given $F \in \widetilde{\mathfrak{F}}$ the invariant
$\mathfrak{K}([\ell^{2,F}])$ satisfies for $(i,j) \in \Delta$
$$\mathfrak{K}([\ell^{2,F}])(i,j)=\Bigg[\prod_{k=i+1}^{j}F(k)
\Bigg].$$ In particular, by noting that there is a well defined
product in $\mathcal{F}$ which is given for $[f_1], [f_2] \in
\mathcal{F}$ by $[f_1] \cdot [f_2]=[f_1 \cdot f_2]$ we obtain the
relations
$$\mathfrak{K}([\ell^{2,F}])(i,j)=\prod_{k=i}^{j-1}
\mathfrak{K}([\ell^{2,F}])(k,k+1).$$ We define an embedding
$$\iota \colon \mathfrak{F} \to \mathrm{Map}(\Delta,\mathcal{F})$$
which is given for $F \in \mathfrak{F}$ by the formula
$$\iota(F)(i,j)=\prod_{k=i+1}^{j}F(k), \quad
(i,j) \in \Delta.$$ By the same formula we define also an embedding
$$\iota_n \colon \mathfrak{F}_n \to
\mathrm{Map}(\Delta_n,\mathcal{F}).$$ As a Corollary of
Corollary~\ref{a2} we obtain the following statement.
\begin{cor}\label{a4}
Every $A \in \iota(\mathfrak{F}) \subset
\mathrm{Map}(\Delta,\mathcal{F})$ is realizable as the invariant of
a scale Hilbert space, i.e.~there exists a scale Hilbert space
$\mathcal{H}$ such that
$$\mathfrak{K}([\mathcal{H}])=A.$$
Similarly, every $A \in \iota_n(\mathfrak{F}_n)$ is realizable as
the invariant of a scale Hilbert $n$-tuple.
\end{cor}
Our second main result deals with the question if there are other
invariants than $\iota(\mathfrak{F})$ which can be realized by scale
Hilbert spaces. It actually deals with the question of scale Hilbert
triples, but these can be used to construct new scale Hilbert
spaces. For triples the set $\Delta_3$ has just cardinality three,
namely
$$\Delta_3=\{(0,1),(1,2),(0,2)\}.$$
We identify $\mathrm{Map}(\Delta_3,\mathcal{F})$ with
$\mathcal{F}^3$ via the identification
$$\mathrm{Map}(\Delta_3,\mathcal{F}) \to \mathcal{F}^3,
\quad F \mapsto \big(F(0,1),F(1,2),F(0,2)\big).$$ For given $\phi_1,
\phi_2 \in \mathcal{F}$ we introduce the set
$$\mathcal{B}(\phi_1,\phi_2)=\big\{\phi_3 \in \mathcal{F}:
\exists\,\,H \in \mathscr{S}_3,\,\,
\mathfrak{K}(H)=(\phi_1,\phi_2,\phi_3)\big\}.$$ We know from
Corollary~\ref{a4} that
$$\phi_1 \cdot \phi_2 \in \mathcal{B}(\phi_1,\phi_2)$$
so that in particular, $\mathcal{B}(\phi_1,\phi_2)$ is not empty.
Our second main result is the following Theorem.
\\ \\
\textbf{Theorem~B } \emph{For any $\phi_1,\phi_2 \in \mathcal{F}$,
the set $\mathcal{B}(\phi_1,\phi_2)$ has infinite cardinality.}
\\ \\
As a Corollary of Theorem~B we get the following assertion which we
already discussed in Remark~\ref{nisu}.
\begin{cor}\label{b1}
The map $\mathfrak{J} \colon \mathfrak{F} \to \mathscr{S}$  is not
surjective, as neither the maps $\mathfrak{J}_n \colon
\mathfrak{F}_n \to \mathscr{S}_n$ for every $n \geq 3$.
\end{cor}

\section[Proof of Theorem~A and its Corollaries]
{Proof of Theorem~A and its Corollaries}

The proof of Theorem~A is based on three Lemmas which we prove
first.

\begin{lemma}\label{exi}
If $f \in \widetilde{\mathcal{F}}$, then the tuple
$(\ell^2,\ell^2_f)$ is a scale Hilbert pair.
\end{lemma}
\textbf{Proof: } Abbreviate by $I \colon \ell^2_f \to \ell^2$ the
inclusion. We first observe that the inclusion is continuous.
Indeed, by the assumption that $f$ is monotone we obtain
$$||Ix||_{\ell^2} \leq \frac{1}{\sqrt{f(1)}}
||x||_{\ell^2_f}.$$ To show that $I$ is compact, we denote for $n
\in \N$ by
$$ \Pi_n \colon \ell^2 \to \ell^2$$ the orthogonal projection of a
sequence to its first $n$ entries,
$$(x_1, \ldots,
x_n,x_{n+1},\ldots) \mapsto (x_1,\ldots x_n, 0,\ldots).$$ The
operators
$$I_n =\Pi_n \circ I \colon \ell^2_f \to \ell^2$$
have finite dimensional image and are therefore compact. Since $f$
is monotone we obtain
$$||I-I_n||_{\mathcal{L}(\ell^2_f,\ell^2)} = \frac{1}{\sqrt{f(n+1)}}$$
where $||\cdot||_{\mathcal{L}}$ denotes the operator norm. Because
$f$ is unbounded, we conclude that $I$ is the uniform limit of
compact operators and therefore itself compact. This proves
condition (i) in the Definition of a scale Hilbert pair.
\\
It remains to check condition (ii) in the Definition of a scale
Hilbert pair, i.e.~that $\ell^2_f$ is dense in $\ell^2$. To see that
let $x \in \ell^2$ and define $x^n=\Pi_n x$ for $n \in \N$. We note
that $x^n \in \ell^2_f$ and the sequence $x^n$ converges to $x$ in
the $\ell^2$-norm as $n$ goes to infinity. This proves (ii) and
hence the Lemma. \hfill $\square$

\begin{lemma}\label{uniq}
Assume that $f_1, f_2 \in \widetilde{\mathcal{F}}$. Then the scale
Hilbert pairs $(\ell^2,\ell^2_{f_1})$ and $(\ell^2,\ell^2_{f_2})$
are scale isomorphic iff $f_1 \sim f_2$.
\end{lemma}
\textbf{Proof: } We first prove the implication $"\Rightarrow"$.
Assume that $(\ell^2,\ell^2_{f_1})$ and $(\ell^2,\ell^2_{f_2})$ are
scale isomorphic. Then there exists a scale isomorphism
$$\Phi \colon (\ell^2,\ell^2_{f_1}) \to (\ell^2,\ell^2_{f_2})$$
with inverse
$$\Psi \colon (\ell^2,\ell^2_{f_2}) \to (\ell^2, \ell^2_{f_1}).$$
We abbreviate
$$c=\max\Big\{||\Phi||_{\mathcal{L}(\ell^2,\ell^2)},
||\Phi||_{\mathcal{L}(\ell^2_{f_1},\ell^2_{f_2})},
||\Psi||_{\mathcal{L}(\ell^2,\ell^2)},
||\Psi||_{\mathcal{L}(\ell^2_{f_2},\ell^2_{f_1})}\Big\}$$ where
$||\cdot||_{\mathcal{L}}$ is the operator norm. As in the proof of
Lemma~\ref{exi} we denote for $n \in \N$ by
$$ \Pi_n\colon \ell^2 \to \ell^2$$
the orthogonal projection of a sequence to its first $n$ entries,
$$(x_1, \ldots,
x_n,x_{n+1},\ldots) \mapsto (x_1,\ldots ,x_n, 0,\ldots).$$ For $n,m
\in \N$ we introduce the map
$$A^n_m \colon \Pi_n \ell^2 \to \Pi_n \ell^2, \qquad
A^n_m =\Pi_n \circ \Psi \circ \Pi_{m-1} \circ \Phi.$$ We first prove
the following Claim.
\\ \\
\textbf{Claim: } \emph{Assume that $n,m \in \N$ satisfy
$f_1(n)<\frac{f_2(m)}{c^4}$, then the map $A^n_m$ is a bijection.}
\\ \\
To prove the Claim assume that $\xi$ is in the kernel of $A^n_m$,
i.e. $\xi \in \ell^2$ satisfies
$$\Pi_n \xi=\xi, \quad A^n_m \xi=0.$$
It follows that
\begin{eqnarray*}
\xi=\Pi_n \xi =\Pi_n \Psi \Phi \xi=\Pi_n \Psi
(\mathrm{id}-\Pi_{m-1}) \Phi \xi.
\end{eqnarray*}
We estimate
\begin{eqnarray*}
||\xi||_{\ell^2}&=&||\Pi_n \Psi (\mathrm{id}-\Pi_{m-1}) \Phi
\xi||_{\ell^2}\\
&\leq& c||(\mathrm{id}-\Pi_{m-1}) \Phi \xi||_{\ell^2}\\
&\leq& \frac{c}{\sqrt{f_2(m)}} ||(\mathrm{id}-\Pi_{m-1}) \Phi \xi||_{\ell^2_{f_2}}\\
&\leq& \frac{c^2}{\sqrt{f_2(m)}}||\xi||_{\ell^2_{f_1}}\\
&\leq& \frac{c^2 \sqrt{f_1(n)}}{\sqrt{f_2(m)}} ||\xi||_{\ell^2}.
\end{eqnarray*}
But by the assumption of the claim
$$\frac{c^2 \sqrt{f_1(n)}}{\sqrt{f_2(m)}}<1$$
implying that
$$||\xi||_{\ell^2}=0$$
and hence
$$\xi=0.$$
This proves that $A^n_m$ is injective and since it is an
endomorphism of a finite dimensional vector space we conclude that
$A^n_m$ is a bijection. This finishes the proof of the Claim.
\\ \\
We next show how the Claim can be used to deduce the implication
$"\Rightarrow"$ of the Lemma. We continue assuming the hypothesis of
the Claim. Since the map $A^n_m$ factors as
$$A^n_m=(\Pi_n \circ \Psi)(\Pi_{m-1} \circ \Phi)$$
and $A^n_m$ is bijective by the Claim, we deduce that $\Pi_n \circ
\Psi|_{\Pi_{m-1}\ell^2}$ is surjective. Hence we obtain
$$n=\mathrm{dim} \big(\Pi_n \ell^2\big)=\mathrm{dim}
\big(\mathrm{im} \Pi_n \Psi|_{\Pi_{m-1}\ell^2}\big) \leq
\mathrm{dim}\big( \Pi_{m-1} \ell^2\big)=m-1<m.$$ Hence we have shown
the implication
$$f_1(n) < \frac{f_2(m)}{c^4} \quad \Longrightarrow \quad n<m.$$
We conclude from this that the inequality
$$f_1(n) \geq \frac{f_2(n)}{c^4}$$
has to hold for each $n \in \N$. Interchanging the roles of $\Psi$
and $\Phi$ we obtain the reverse inequality
$$f_2(n) \geq \frac{f_1(n)}{c^4}.$$
This proves that $f_1$ and $f_2$ are equivalent.
\\ \\
We are left with showing the inverse implication $"\Leftarrow"$ of
the Lemma. Hence we assume that $f_1 \sim f_2$. But under this
assumption $\mathrm{id}|_{\ell^2}$ gives a scale isomorphism between
$(\ell^2,\ell^2_{f_1})$ and $(\ell^2,\ell^2_{f_2})$. Hence
$(\ell^2,\ell^2_{f_1})$ and $(\ell^2,\ell^2_{f_2})$ are scale
isomorphic. This finishes the proof of the Lemma. \hfill $\square$ \
\\ \\
Before stating the next Lemma we first introduce another equivalence
relation for scale Hilbert spaces different from scale isomorphism.

\begin{fed}
A \emph{scale isometry} $\Phi$ from a scale Hilbert space
$\mathcal{H}=\{(H_k,\langle \cdot,\cdot \rangle_k)\}$ to a scale
Hilbert space $\mathcal{H}'=\{(H'_k,\langle \cdot, \cdot
\rangle_k)\}$ is a linear map $\Phi \colon H_0 \to H_0'$ which
restricts for all $k \in \N_0$ to an isometry $\Phi|_k \colon H_k
\to H'_k$. Two scale Hilbert spaces are called \emph{scale
isometric}, if there exists a scale isometry between them.
\end{fed}
Note that a scale isometry is a special case of a scale isomorphism,
so that two scale isometric scale Hilbert spaces are also scale
isomorphic. Moreover, the same definition also applies to scale
Hilbert $n$-tuples for any $n \in \N$.

\begin{lemma}\label{sur}
Let $(H,W)$ be an infinite dimensional scale Hilbert pair. Then
there exists a unique $f \in \widetilde{\mathcal{F}}$ such that
$(H,W)$ is scale isometric to $(\ell^2,\ell^2_f)$.
\end{lemma}
\textbf{Proof: } By the Riesz representation theorem there exists a
bounded linear operator $A \colon W \to W$ such that
$$\langle w_1, w_2 \rangle_H=\langle w_1, A w_2 \rangle_W, \quad
w_1, w_2 \in W.$$ The operator $A$ is symmetric and we next show
that it is compact. Choose a sequence $w_\nu$ in the unit ball of
$W$, i.e.
$$||w_\nu||_W \leq 1 , \quad \nu \in \N.$$
Since the inclusion $W \hookrightarrow H$ is compact we deduce that
$w_\nu$ has a convergent subsequence $w_{\nu_j}$ in $H$. In
particular, $w_{\nu_j}$ is a Cauchy sequence in $H$. We claim that
$A w_{\nu_j}$ is a Cauchy sequence in $W$. Denote by $||A||>0$ the
operator norm of the bounded linear operator $A \colon W \to W$.
Since $w_{\nu_j}$ is a Cauchy sequence in $W$ there exists for given
$\epsilon>0$ a positive integer $j_0=j_0(\epsilon) \in \N$ such that
$$||w_{\nu_j}-w_{\nu_{j'}}||_H \leq \frac{\epsilon}{\sqrt{||A||}},
\quad j,j' \geq j_0.$$
We further abbreviate
$$v=w_{\nu_j}-w_{\nu_{j'}}.$$
We estimate
\begin{eqnarray*}
0 &\leq& \bigg\langle v-\frac{1}{||A||}Av,v-\frac{1}{||A||} Av
\bigg\rangle_H\\
&=&||v||^2_H-\frac{2}{||A||}\langle Av,v \rangle_H
+\frac{1}{||A||^2} \langle Av, Av \rangle_H\\
&=&||v||^2_H-\frac{2}{||A||}\langle Av, Av
\rangle_W+\frac{1}{||A||^2} \langle Av, A^2 v \rangle_W\\
&\leq& \frac{\epsilon^2}{||A||}-\frac{2}{||A||}||Av||^2_W
+\frac{1}{||A||^2}
||Av||_W||A^2v||_W\\
&\leq&
\frac{\epsilon^2}{||A||}-\frac{2}{||A||}||Av||^2_W+\frac{1}{||A||}||Av||^2_W\\
&=&\frac{\epsilon^2}{||A||}-\frac{1}{||A||}||Av||^2_W
\end{eqnarray*}
from which we conclude
$$||Aw_{\nu_j}-Aw_{\nu_{j'}}||_W=||Av||_W \leq \epsilon.$$
This proves that $Aw_{\nu_j}$ is a Cauchy sequence in $W$ and since
$W$ is complete it has to converge. We deduce that $A$ is a compact
operator.

We next apply the spectral theorem to the compact symmetric operator
$A$. Since $A$ is further positive we conclude that there exists an
orthogonal Schauder basis $\{e_n\}_{n \in \N}$ of $W$ with the
following properties.
\begin{description}
 \item[(i)] For each $n \in \N$ the vector $e_n$ is an eigenvector
  of $A$ to a real eigenvalue $\lambda_n>0$.
 \item[(ii)] The eigenvalues $\lambda_n$ build a monotone decreasing
 sequence.
\end{description}
Since $\{e_n\}_{n \in \N}$ is an orthogonal Schauder basis of $W$ we
can represent each $w \in W$ in the form
\begin{equation}\label{W}
w=\sum_{n=1}^\infty x_n e_n, \quad x=(x_1,x_2,\cdots) \in \ell^2.
\end{equation}
We next construct an orthogonal basis for $H$. For $n \in \N$ define
$$\bar{e}_n:=\frac{1}{\sqrt{\lambda_n}} e_n.$$
Denoting for $n,m \in \N$ by $\delta^m_n$ the Kronecker symbol we
compute
$$\langle \bar{e}_n,\bar{e}_m \rangle_H=
\langle \bar{e}_n, A \bar{e}_m \rangle_W=
\frac{\lambda_m}{\sqrt{\lambda_m \lambda_n}}\langle e_n, e_m
\rangle_W=\frac{\lambda_m}{\sqrt{\lambda_m \lambda_n}}\delta_n^m
=\delta_n^m$$ and hence the vectors $\bar{e}_n$ are orthogonal to
each other. To see that they form a Schauder basis of $H$ define
$$H'=\overline{\mathrm{span}(\bar{e}_1,\bar{e}_2,\cdots)}^{||\cdot||_H}$$
to be the $||\cdot||_H$-closure of the vector space spanned by
$\{\bar{e}_n\}_{n \in \N}$. We observe that $H'$ is a closed
subspace of $H$ and $W$ is dense in $H'$. Since $W$ is dense in $H$
by assumption we conclude that
$$H'=H.$$
We now define an isometry
$$\Phi \colon H \to \ell^2$$
in the following way. By the reasoning above each element $h \in H$
has a unique representation
$$h=\sum_{n=1}^\infty y_n \bar{e}_n, \quad
y=(y_1,y_2,\ldots) \in \ell^2$$ and we set
$$\Phi(h)=y.$$
We next study the restriction of $\Phi$ to $W$. Define $f \in
\widetilde{\mathcal{F}}$ by
$$f(n)=\frac{1}{\lambda_n}, \quad n \in \N.$$
Since $\lambda_n$ is a monotone decreasing zero sequence, the
function $f$ is actually monotone and unbounded. We claim that the
restriction of $\Phi$ to $W$ gives an isometry
$$\Phi|_W \colon W \to \ell^2_f \subset \ell^2.$$ To prove that
assertion let $w^1, w^2 \in W$. By (\ref{W}) there exist
$x^1=(x^1_1, x^1_2, \cdots) \in \ell^2$ and $x^2=(x^2_1, x^2_2,
\cdots) \in \ell^2$ such that for $i \in \{1,2\}$ we have
$$w^i=\sum_{n=1}^\infty x^i_n e_n=\sum_{n=1}^\infty
\sqrt{\lambda_n}x^i_n \bar{e}_n.$$ In particular, we get
$$\Phi(w^i)=(\sqrt{\lambda}_1 x^i_1, \sqrt{\lambda}_2 x^i_2,
\cdots).$$ Hence we compute
$$\langle \Phi(w^1), \Phi(w^2)\rangle_{\ell^2_f}=
\sum_{n=1}^\infty f(n) \lambda_n x^1_n x^2_n =\sum_{n=1}^\infty
x^1_n x^2_n =\langle w^1, w^2 \rangle_W.$$ This proves that
$\Phi|_W$ interchanges the two scalar products. In particular,
$\Phi|_W$ is injective. To see that it is surjective we note that if
$y=(y_1, y_2, \cdots) \in \ell^2_f$, then
$$\bigg(\frac{y_1}{\sqrt{\lambda_1}},\frac{y_2}{\sqrt{\lambda_2}},\cdots
\bigg) \in \ell^2$$ and hence
$$w=\sum_{n=1}^\infty \frac{y_n}{\sqrt{\lambda_n}} e_n \in W.$$
But
$$\Phi(w)=y$$
which shows that $\Phi|_{W}:W \to \ell^2_f$ is surjective. This
finishes the proof that $\Phi|_{W}$ is an isometry from $W$ to
$\ell_f^2$. In particular,
$$\Phi \colon (H,W) \to (\ell^2,\ell^2_f)$$
defines a scale isometry.

It finally remains to show that $f \in \widetilde{\mathcal{F}}$ is
unique with this property. To see this assume that $f_1, f_2 \in
\widetilde{\mathcal{F}}$ such that there exist scale isometries
$$\Phi_1 \colon (H,W) \to (\ell^2,\ell^2_{f_1}), \quad
\Phi_2 \colon (H,W) \to (\ell^2,\ell^2_{f_2}).$$ Then
$$\Psi=\Phi_2 \circ \Phi_1^{-1} \colon (\ell^2,\ell^2_{f_1})
\to (\ell^2,\ell^2_{f_2})$$ is also a scale isometry. Let
$\{\varepsilon_n\}_{n\in \N}$ be the standard basis of $\ell^2$
given by
$$\varepsilon_n=(\delta_n^1, \delta_n^2, \cdots).$$
For $f \in \widetilde{\mathcal{F}}$ and $n \in \N$ we set
$$\varepsilon^f_n=\frac{1}{\sqrt{f(n)}}\varepsilon_n$$
for the standard $\ell^2$-basis of $\ell^2_f$. We further define by
$$A^f \colon \ell^2_f \to \ell^2_f$$
the linear map which is given on basis vectors by
\begin{equation}\label{af}
A^f \varepsilon^f_n=\frac{1}{f(n)}\varepsilon^f_n.
\end{equation}
With this convention we have for vectors $w_1, w_2 \in \ell^2_f$ the
equality
$$\langle w_1, w_2 \rangle_{\ell^2}=\langle w_1, A^f w_2
\rangle_{\ell^2_f}.$$ Now using that $\Psi$ is a scale isometry we
compute for $w_1, w_2 \in \ell^2_{f_1}$
\begin{eqnarray*}
\big\langle \Psi w_1, A^{f_2} \Psi w_2 \big\rangle_{\ell^2_{f_2}}
&=&\big\langle \Psi w_1, \Psi w_2 \big\rangle_{\ell^2}\\
&=&\big\langle w_1,  w_2 \big\rangle_{\ell^2}\\
&=&\big\langle w_1, A^{f_1}w_2 \big\rangle_{\ell^2_{f_1}}\\
&=&\big\langle \Psi w_1, \Psi A^{f_1}w_2 \big\rangle_{\ell^2_{f_2}}
\end{eqnarray*}
implying that
$$A^{f_2}\Psi=\Psi A^{f_1}.$$
This shows that $A^{f_1}$ and $A^{f_2}$ have the same eigenvalues
and by (\ref{af}) we deduce the following equality of sets
$$\{f_1(n): n \in \N\}=\{f_2(n): n \in \N\}.$$
Since $f_1$ and $f_2$ are monotone we get
$$f_1(n)=f_2(n), \quad n \in \N.$$
This proves the uniqueness part and hence the Lemma follows. \hfill
$\square$
\\ \\
\textbf{Proof of Theorem~A: }By Lemma ~\ref{exi} the map
$$\widehat{\mathfrak{J}}_2 \colon \widetilde{\mathcal{F}}
\to \mathscr{S}_2, \quad f \mapsto [(\ell^2,\ell^2_f)]$$ is well
defined. By Lemma~\ref{uniq} this map induces a map
$$\mathfrak{J}_2 \colon \mathcal{F} \to \mathscr{S}_2.$$
By Lemma~\ref{sur} the map $\mathfrak{J}_2$ is surjective and again
by Lemma~\ref{uniq} it is also injective. This finishes the proof of
Theorem~A. \hfill $\square$
\\ \\
\textbf{Proof of Corollary~\ref{a1}: }Note that $\{(H_k,\langle
\cdot,\cdot\rangle_k),(H_{k-1},\langle \cdot,\cdot \rangle_{k-1})\}$
is a scale Hilbert pair. Hence Theorem~A implies the Corollary.
\hfill $\square$
\\ \\
\textbf{Proof of Corollary~\ref{a3}: }Assume that we have given two
infinite dimensional scale Hilbert spaces
$\mathcal{H}=\{(H_k,\langle \cdot,\cdot \rangle_k)\}$ and
$\mathcal{H}'=\{(H_k',\langle \cdot,\cdot \rangle'_k)\}$ which are
scale isomorphic to each other. In particular, there exists a scale
isomorphism
$$\Phi \colon \mathcal{H} \to \mathcal{H}'.$$
If $(i,j) \in \Delta$ then by restricting $\Phi$ we obtain a scale
isomorphism for scale Hilbert pairs
$$\Phi_{i,j} \colon \big((H_i,\langle \cdot,\cdot \rangle_i),
(H_j,\langle \cdot, \cdot \rangle_j)\big) \to \big((H_i',\langle
\cdot,\cdot \rangle'_i), (H_j',\langle \cdot, \cdot
\rangle'_j)\big).$$ In particular, we conclude
$$\big[(H_i,\langle \cdot,\cdot \rangle_i),
(H_j,\langle \cdot, \cdot \rangle_j)\big] = \big[(H_i',\langle
\cdot,\cdot \rangle'_i), (H_j',\langle \cdot, \cdot \rangle'_j)\big]
\in \mathscr{S}_2.$$ Hence the map $\mathfrak{K}$ is well defined,
since the map $\mathfrak{J}_2$ is well defined by Theorem~A. The
same reasoning also applies to $\mathfrak{K}_n$ for $n$ an integer
greater than $1$. This finishes the proof of the Corollary. \hfill
$\square$
\\ \\
\textbf{Proof of Corollary~\ref{a2}: } For the proof of the
Corollary we use the convention
$$\mathscr{S}_\infty=\mathscr{S},\quad
\mathfrak{F}_\infty=\mathfrak{F}, \quad
\mathfrak{J}_\infty=\mathfrak{J}, \quad \mathfrak{K}_\infty
=\mathfrak{K}$$ and we assume that $n \in \{2,3,\ldots,\infty\}$. We
first show that the map $\mathfrak{J}_n$ is well defined, i.e.~that
$[\ell^{2,F}] \in \mathscr{S}_n$ for every $F \in
\widetilde{\mathfrak{F}}_n$. We claim that if $k$ is a positive
integer less than $n$ that the inclusion $\ell_k^{2,F}
\hookrightarrow \ell_{k-1}^{2,F}$ is compact. If $k=1$ this
inclusion corresponds to the inclusion $\ell_{F(1)}^2
\hookrightarrow \ell^2$ which is compact by Lemma~\ref{exi}. Now
assume that $k>1$. Let $\{\varepsilon_\nu\}_{\nu \in \N}$ be the
standard orthogonal basis of $\ell^2$. Let
$\{\widetilde{\varepsilon}_{\nu}\}_{\nu \in \N}$ the orthogonal
basis of $\ell^{2,F}_{k-1}=\ell^2_{\prod_{j=1}^{k-1} F(j)}$ defined
by
$$\widetilde{\varepsilon}_\nu=\frac{1}{\sqrt{\prod_{j=1}^{k-1}F(j)(\nu)}}
\varepsilon_\nu.$$ Denote by
$$I \colon \ell^{2,F}_{k-1} \to \ell^2$$
the isometry which is given on basis vectors by
$$I(\widetilde{\varepsilon}_\nu)=\varepsilon_\nu, \quad
\nu \in \N.$$ Note that the restriction of $I$ to $\ell^{2,F}_k$
gives an isometry
$$I|_{\ell^{2,F}_k}: \ell^{2,F}_k \to \ell^2_{F(k)}.$$
Hence we conclude that the pair $(\ell^{2,F}_{k-1},\ell^{2,F}_k)$ is
scale isometric to the pair $(\ell^2,\ell^2_{F(k)})$ and the
compactness of the embedding follows again from Lemma~\ref{exi}.
\\
We next show that the intersection $\bigcap_{j=0}^{n-1}\ell^{2,F}_j$
is dense in $\ell^{2,F}_k$ for every nonnegative integer $k$ less
than $n$. To see that let
$$\mathfrak{f}=\mathrm{span}\{\varepsilon_\nu: \nu \in \N\}$$
be the subspace of $\ell^2$ consisting of finite linear combinations
of the standard basis vectors of $\ell^2$. We note that
$$\mathfrak{f} \subset \ell^{2,F}_k$$
is a dense subspace for every nonnegative integer $k$ less than $n$.
In particular,
$$\mathfrak{f} \subset \bigcap_{j=0}^{n-1}\ell^{2,F}_j.$$
This shows that $\bigcap_{j=0}^{n-1}\ell^{2,F}_j$ is dense in
$\ell^{2,F}_k$. We conclude that $[\ell^{2,F}] \in \mathscr{S}_n$
and hence the map $\mathfrak{J}_n$ is well defined.
\\
We are left with showing injectivity of the map $\mathfrak{J}_n$.
Hence assume that $F_1,F_2 \in \widetilde{\mathfrak{F}}_n$ such that
$$[F_1] \neq [F_2] \in \mathfrak{F}.$$
This implies that there exists a positive
integer $k$ less than $n$ such that
$$[F_1(k)] \neq [F_2(k)] \in \mathcal{F}.$$
We noted before that the pairs $(\ell^{2,F_i}_{k-1},\ell^{2,F_i}_k)$
and $(\ell^2,\ell^2_{F_i(k)})$ are scale isometric for $i \in
\{1,2\}$. In particular, these pairs are scale isomorphic so that we
obtain
$$\mathfrak{J}_2^{-1}\big(\ell^{2,F_i}_{k-1},\ell^{2,F_i}_k\big)
=[F_i(k)], \quad i \in \{1,2\}.$$ Combining the above two facts we
conclude
$$\mathfrak{K}([\ell^{2,F_1}])(k-1,k)=
[F_1(k)] \neq [F_2(k)]= \mathfrak{K}([\ell^{2,F_2}])(k-1,k).$$ In
particular,
$$\mathfrak{K}([\ell^{2,F_1}]) \neq \mathfrak{K}([\ell^{2,F_2}])$$
implying that
$$[\ell^{2,F_1}] \neq [\ell^{2,F_2}]$$
which proves that $\mathfrak{J}_n$ is injective. This finishes the
proof of the Corollary. \hfill $\square$
\\ \\
\textbf{Proof of Corollary~\ref{a4}: }If $A \in
\iota(\mathfrak{F})$, then there exists $F \in
\widetilde{\mathfrak{F}}$ such that
$$A=\iota([F]).$$
Set
$$\mathcal{H}=\ell^{2,F}.$$
Then
$$\mathfrak{K}([\mathcal{H}])=A.$$
The same reasoning applies to scale Hilbert $n$-tuples. \hfill
$\square$

\begin{rem}
\emph{The uniqueness statement in Lemma~\ref{sur} was actually not
used in the proof of Theorem~A. However, it can be used to describe
the set of scale Hilbert pairs modulo the equivalence given by scale
isometry instead of scale isomorphism. Introduce the set
$$\widetilde{\mathscr{S}}_2=\big\{\mathcal{H}\,\,
\textrm{scale Hilbert pair},\,\dim(H_0)=\infty\}/\sim'$$ where
$\sim'$ is the equivalence relation given by scale isometry. Then
the map
$$\widetilde{\mathfrak{J}}_2 \colon \widetilde{\mathcal{F}}
\to \widetilde{\mathscr{S}}_2, \quad f \mapsto (\ell^2,\ell^2_f)$$
gives a bijection between $\widetilde{\mathscr{S}}_2$ and
$\widetilde{\mathcal{F}}$.}
\end{rem}

\section[Proof of Theorem~B and its Corollary]{Proof of Theorem~B
and its Corollary}

The proof of Theorem~B roots on the following idea. Choose
representatives $f_1,f_2 \in \widetilde{\mathcal{F}}$ for $\phi_1$
respectively $\phi_2$. The separable Hilbert space $\ell^2_{f_1}$ is
isometric to $\ell^2$ and we have a canonical isometry
$$I \colon \ell^2_{f_1} \to \ell^2.$$
Let
$$\Phi \colon \ell^2 \to \ell^2$$
be an isometry of $\ell^2$ to itself. Now consider the scale Hilbert
triple
$$\mathcal{H}=\big(\ell^2,\ell^2_{f_1}, I^{-1} \Phi(\ell^2_{f_2})\big).$$
Applying $\Phi^{-1} \circ I$ to $\ell^2_{f_1}$ we get a scale
isomorphism
$$\Phi^{-1} \circ I \colon \big(\ell^2_{f_1}, I^{-1}
\Phi(\ell^2_{f_2})\big) \to \big(\ell^2,\ell^2_{f_2}\big).$$ In
particular,
$$\mathfrak{K}([\mathcal{H}])(1,2)=[f_2]=\phi_2.$$
Moreover, we have
$$\mathfrak{K}([\mathcal{H}])(0,1)=[f_1]=\phi_1.$$
On the other hand $\mathfrak{K}([\mathcal{H}])(0,2)$ depends on
$\Phi$ and we show that by varying $\Phi$ we can achieve infinitely
many values for $\mathfrak{K}([\mathcal{H}])(0,2)$ in the set
$\mathcal{F}$.

We now start with the preparations for the proof of Theorem~B. We
denote by $\mathcal{U}$ the set of all functions $u \colon \N \to
(0,\infty)$ satisfying $\lim_{n \to \infty}u(n)=\infty$. Obviously,
$$\widetilde{\mathcal{F}} \subset \mathcal{U}.$$
We further introduce
$$\mathfrak{S}=\{\sigma \colon \N \to \N:
\sigma\,\,\textrm{bijective}\}$$ the group of permutations of $\N$.
\begin{lemma}
The group $\mathfrak{S}$ acts on $\mathcal{U}$ by
$$\sigma_* u(n)=u(\sigma(n)), \quad \sigma \in \mathfrak{S},\,\, u \in
\mathcal{U},\,\,n \in \N.$$
\end{lemma}
\textbf{Proof: }We prove that the action is well defined, i.e.~that
$\sigma_* u \in \mathcal{U}$. We have to show that
\begin{equation}\label{limu}
\lim_{n \to \infty}\sigma_* u(n)=\infty.
\end{equation}
Pick $r \in \R$. Since $u \in \mathcal{U}$ there exists $n_0=n_0(r)$
such that
\begin{equation}\label{gross}
u(n) \geq r, \quad \forall\,\,n \geq n_0.
\end{equation}
Since $\sigma$ is bijective the set $\{n \in \N: \sigma^{-1}(n) <
n_0\}$ is finite. Hence we can set
$$N_0 :=\max\{n \in \N: \sigma^{-1}(n) < n_0\}.$$
In particular, we have the implication
$$n \geq N_0 \quad \Longrightarrow \quad \sigma(n) \geq n_0.$$
Hence using (\ref{gross}) we conclude
$$\sigma_*u(n)=u(\sigma(n)) \geq r, \quad \forall\,\, n \geq N_0.$$
This proves (\ref{limu}) and hence the Lemma. \hfill $\square$
\begin{lemma}\label{proj}
If $u \in \mathcal{U}$ there exists $\sigma \in \mathfrak{S}$ such
that $\sigma_* u \in \widetilde{\mathcal{F}}$. Moreover, if $\sigma'
\in \mathfrak{S}$ is another element with this property, than
$\sigma_*u=\sigma'_* u$.
\end{lemma}
\begin{rem} Although $\sigma_* u$ in Lemma~\ref{proj} is canonical,
the permutation $\sigma$ need not be. It is only canonical if
$\sigma_* u$ is strictly monotone.
\end{rem}
\textbf{Proof of Lemma~\ref{proj}: }Pick $u \in \mathcal{U}$. We
first note that since $u$ converges to infinity it follows that for
each finite subset $B \subset \N$ the infimum of the restriction of
$u$ to $\N \setminus B$ is attained so that we are allowed to put
$$a_B:=\min\big\{u(n): n \in \N \setminus B\big\}.$$
We set
$$B_0=\emptyset$$
and define recursively for $k \in \N$
$$a_k:=a_{B_{k-1}}, \quad \sigma(k):=\min\big\{n \in \N \setminus B_{k-1}:
u(n)=a_k\big\}, \quad B_k:=B_{k-1} \cup \{\sigma(k)\}.$$ We claim
that
\begin{equation}\label{sig}
\sigma \in \mathfrak{S},
\end{equation}
i.e.~that $\sigma$ is bijective. We first show injectivity. We
assume by contradiction that there exist $k,k' \in \N$ such that
$$\sigma(k)=\sigma(k'), \quad k \neq k'.$$
We can assume without loss of generality that
$$k<k'.$$
It follows from the definition of $B_k$ that
$$B_k=\{\sigma(j): 1 \leq j \leq k\}.$$
We deduce from the definition of $\sigma(k')$ that
$$\sigma(k') \in \N \setminus \{\sigma(j): 1 \leq j \leq k'-1\}
\subset \N \setminus \{\sigma(k)\}=\N \setminus \{\sigma(k')\}$$
which is absurd. Therefore injectivity of $\sigma$ has to hold. We
next show surjectivity again by contradiction. We assume that there
exists $m \in \N$ such that
$$\sigma(k) \neq m, \quad \forall\,\,k \in \N.$$
It follows that
$$m \in \N \setminus B_k,\quad \forall\,\,k \in \N.$$
Therefore
$$a_k \leq u(m), \quad \forall\,\,k \in \N.$$
We conclude
$$u(\sigma(k)) \leq u(m), \quad \forall\,\,k \in \N.$$
Since $\sigma$ is injective as we have already shown we deduce that
$$\#\{n \in \N: u(\sigma(n)) \leq u(m)\}=\infty.$$
But this contradicts the fact that $u$ converges to infinity. Hence
$\sigma$ has to be surjective and (\ref{sig}) is proved.
\\
We next check that $\sigma_* u \in \widetilde{\mathcal{F}}$,
i.e.~$\sigma_* u$ is monotone. To see that we estimate for $k \in
\N$
\begin{eqnarray*}
\sigma_*u(k+1)&=&u(\sigma(k+1))\\
&=&a_{k+1}\\
&=&a_{B_k}\\
&=&\min\big\{u(n): n \in \N \setminus B_k\big\}\\
&=&\min\big\{u(n): n \in \N \setminus \{\sigma(j): 1 \leq j \leq
k\}\big\}\\
&\geq&\min\big\{u(n): n \in \N \setminus \{\sigma(j): 1 \leq j \leq
k-1\}\big\}\\
&=&\sigma_*u(k).
\end{eqnarray*}
This proves monotonicity and hence the existence statement of the
Lemma is settled.
\\
We are left with proving the uniqueness statement of the Lemma. We
prove by induction on $k$ that
\begin{equation}\label{suu}
\sigma_* u(k)=\sigma'_* u(k).
\end{equation}
Using the monotonicity of $\sigma_* u$ and $\sigma'_* u$ and the
bijectivity of $\sigma$ and $\sigma'$ we compute
\begin{eqnarray*}
\sigma_*u(1)&=&\min\{\sigma_* u(n): n \in \N\}\\
&=&\min\{u(\sigma(n)): n \in \N\}\\
&=&\min\{u(n): n \in \N\}\\
&=&\sigma'_*u(1).
\end{eqnarray*}
which is (\ref{suu}) for $k=1$. Assuming (\ref{suu}) for all $j \leq
k$ we obtain
\begin{eqnarray*}
\sigma_*u(k+1)&=&\min\{\sigma_* u(n):n \geq k+1\}\\
&=&\min\Big(\{u(n): n \in \N\}\setminus \{\sigma_* u(j): 1\leq j
\leq k\}\Big)\\
&=&\min\Big(\{u(n): n \in \N\}\setminus \{\sigma'_* u(j): 1\leq j
\leq k\}\Big)\\
&=&\sigma'_*u(k+1).
\end{eqnarray*}
We have proved the induction step and hence (\ref{suu}) follows for
all $k \in \N$. This finishes the proof of uniqueness and hence of
the Lemma. \hfill$\square$
\\ \\
By the previous Lemma we obtain a well defined map
$$P \colon \mathcal{U} \to \widetilde{\mathcal{F}}.$$
Namely, let $u \in \mathcal{U}$ and choose $\sigma \in \mathfrak{S}$
such that $\sigma_* u \in \widetilde{\mathcal{F}}$ and set
$$P(u)=\sigma_* u.$$
The uniqueness statement of the Lemma assures that $P$ is well
defined, i.e.~independent of the choice of $\sigma$. Moreover, we
have the following Corollary.
\begin{cor} The map $P \colon \mathcal{U} \to
\widetilde{\mathcal{F}}$ is a projection, i.e. $P^2=P$.
\end{cor}
\textbf{Proof: } Since $Pu \in \widetilde{\mathcal{F}}$, we have
$(\mathrm{id})_* (Pu) \in \widetilde{\mathcal{F}}$ and hence
$$P^2u=P(Pu)=Pu.$$
This proves the Corollary. \hfill $\square$
\\ \\
For $u_1, u_2 \in \mathcal{U}$ the product is defined pointwise by
$$(u_1 \cdot u_2)(n)=u_1(n) \cdot u_2(n), \quad n \in \N.$$
Note that $u_1 \cdot u_2 \in \mathcal{U}$. For $\sigma \in
\mathfrak{S}$ we define a map
$$\wp_\sigma \colon \widetilde{\mathcal{F}}\times
\widetilde{\mathcal{F}} \to \widetilde{\mathcal{F}}$$ by
$$\wp_\sigma(f_1,f_2)=P(f_1 \cdot \sigma_* f_2), \quad
f_1, f_2 \in \widetilde{\mathcal{F}}.$$ Note that
$$\wp_{\mathrm{id}}(f_1,f_2)=f_1 \cdot f_2, \quad f_1, f_2
\in \widetilde{\mathcal{F}}.$$ If $f \in \widetilde{\mathcal{F}}$ we
denote by $[f]$ the equivalence class of $f$ in $\mathcal{F}$.
\begin{fed}
Given $f_1, f_2 \in \widetilde{\mathcal{F}}$, a subset
$\mathfrak{S}_0 \subset \mathfrak{S}$ is called
\emph{$(f_1,f_2)$-wild}, if
$$[\wp_\sigma(f_1,f_2)] \neq [\wp_{\sigma'}(f_1,f_2)],
  \quad \forall\,\,\sigma, \sigma' \in \mathfrak{S}_0,\,\,\sigma \neq
  \sigma'.$$
\end{fed}
\begin{prop}\label{wilder}
Given $f_1, f_2 \in \widetilde{\mathcal{F}}$, there exists an
$(f_1,f_2)$-wild subset $\mathfrak{S}_0 \subset \mathfrak{S}$ of
infinite cardinality.
\end{prop}
We prove the Proposition with the help of the following Lemma.
\begin{lemma}\label{willi}
Given $f_1, f_2 \in \widetilde{\mathcal{F}}$ and a finite
$(f_1,f_2)$-wild subset $\mathfrak{S}_0 \subset \mathfrak{S}$, then
there exists $\sigma \in \mathfrak{S}\setminus \mathfrak{S}_0$ such
that $\mathfrak{S}_0 \cup \{\sigma\}$ is still an $(f_1,f_2)$-wild
subset of $\mathfrak{S}$.
\end{lemma}
\textbf{Proof: } We prove the Lemma in six steps.
\\ \\
\textbf{Step~1: } \emph{We can assume without loss of generality
that $\mathfrak{S}_0$ is nonempty.}
\\ \\
This follows since $\{\mathrm{id}\}$ is an $(f_1,f_2)$-wild subset
of $\mathfrak{S}$.
\\ \\
\textbf{Step~2: } \emph{The function $g=g_{f_1,f_2} \colon \N \to
(0,\infty)$ which is defined for $n \in \N$ by the formula
$$g_{f_1,f_2}(n)=g(n)=\min\big\{f_1(k) f_2(n+1-k): 1 \leq k \leq n\big\}$$
lies in $\widetilde{\mathcal{F}}$, i.e. $g$ is monotone and
unbounded.}
\\ \\
To prove Step~2 we first show that $g$ is monotone. Let $k \in \{1,
\ldots, n+1\}$ such that
$$g(n+1)=f_1(k)f_2(n+2-k).$$
We first treat the case where $k \leq n$. In this case we estimate
using the monotonicity of $f_2$
$$g(n) \leq f_1(k) f_2(n+1-k) \leq f_1(k)f_2(n+2-k)=g(n+1).$$
If $k=n+1$ we estimate using the monotonicity of $f_1$
$$g(n) \leq f_1(n) f_2(1) \leq f_1(n+1) f_2(1)=g(n+1).$$
We have shown that $g$ is monotone. We next show that $g$ is
unbounded. Since $f_1$ and $f_2$ are unbounded there exists for
given $r \in \R$ a positive integer $n_0=n_0(r)$ with the property
that
$$f_1(n) \geq \frac{r}{\min\{f_1(1),f_2(1)\}}, \quad
f_2(n)\geq \frac{r}{\min\{f_1(1),f_2(1)\}}, \qquad \forall\,\,n \geq
n_0.$$ Using the above inequality and the monotonicity of $f_1$ and
$f_2$ we estimate for $k \in \{1, \ldots,n_0\}$
\begin{equation}\label{small}
f_1(k)f_2(2n_0+1-k) \geq f_1(1) f_2(n_0)\geq r.
\end{equation}
Similarly, we estimate for $k \in \{n_0+1, \ldots,2n_0\}$
\begin{equation}\label{large}
f_1(k)f_2(2n_0+1-k) \geq f_1(n_0) f_2(1)\geq r.
\end{equation}
Inequalities (\ref{small}) and (\ref{large}) imply
$$g(2n_0) \geq r.$$
This proves that $g$ is unbounded and hence Step~2 follows.
\\ \\
\textbf{Step~3: } \emph{Definition of $\sigma \in \mathfrak{S}$.}
\\ \\
For $\ell \in \N_0$ we introduce the shift map
$$s_\ell \colon \widetilde{\mathcal{F}} \to
\widetilde{\mathcal{F}}$$ which is given for $f \in
\widetilde{\mathcal{F}}$ by the formula
$$s_\ell(f)(n)=f(n+\ell), \quad n \in \N.$$
Note that $s_\ell$ is well defined, i.e. $s_\ell(f)$ is still
monotone and unbounded. Since $\mathfrak{S}_0$ is finite and
nonempty by Step~1, we can set for $\ell \in \N$
$$b_\ell=
\max\big\{\wp_{\sigma}(f_1,f_2)(\ell): \sigma \in
\mathfrak{S}_0\big\}.$$ Again for $\ell \in \N$ we further introduce
the set
$$A_\ell=\big\{n \in \N: g_{s_{\ell-1}(f_1),s_{\ell-1}(f_2)}(n) \geq
\ell b_\ell\big\}.$$ Applying Step~2 to
$g_{s_{\ell-1}(f_1),s_{\ell-1}(f_2)}$ we conclude that the set
$A_\ell$ is nonempty. Hence we can set
$$a_\ell=\min\{n: n \in A_\ell\}.$$
We put
$$\ell_1=1$$
and define recursively for $\nu \in \N$
$$
\ell_{\nu+1}=a_{\ell_\nu}+\ell_\nu.$$ Note that for any $\nu \in \N$
$$\ell_\nu <\ell_{\nu+1}.$$
We define $\sigma$ by the formula
$$\sigma(k)=\ell_\nu+\ell_{\nu+1}-k-1, \quad
\ell_\nu \leq k \leq \ell_{\nu+1}-1,\,\,\nu \in \N\\.$$
To show that
$\sigma \in \mathfrak{S}$ we have to check that $\sigma$ is a
bijective map from $\N$ to $\N$. But
$$\sigma|_{\{\ell_\nu, \ldots, \ell_{\nu+1}-1\}} \colon \{\ell_\nu, \ldots,
\ell_{\nu+1}-1\} \to \{\ell_\nu,\ldots, \ell_{\nu+1}-1\}$$ are
bijections for every $\nu \in \N$. This proves that $\sigma$ is a
bijection and finishes Step~3.
\\ \\
\textbf{Step~4: } \emph{For every $\nu \in \N$ we have the
inequality}
\begin{equation}\label{sigi}
(f_1 \cdot \sigma_* f_2)(k) \geq \ell_\nu b_{\ell_\nu}, \quad
\forall\,\,k \geq \ell_\nu.
\end{equation}
We first consider the case where $k \in \{\ell_\nu, \ldots,
\ell_{\nu+1}-1\}$ and estimate
\begin{eqnarray*}
(f_1 \cdot \sigma_* f_2)(k) &=& f_1(k)f_2(\sigma(k))\\
&=&f_1(k)f_2(\ell_\nu+\ell_{\nu+1}-k-1)\\
&=&\big(s_{\ell_\nu-1}(f_1)(k-\ell_\nu+1)\big)
\big(s_{\ell_\nu-1}(f_2)(\ell_{\nu+1}-k)\big)\\
&\geq&g_{s_{\ell_\nu-1}(f_1),s_{\ell_\nu-1}(f_2)}(\ell_{\nu+1}-\ell_\nu)\\
&=&g_{s_{\ell_\nu-1}(f_1),s_{\ell_\nu-1}(f_2)}(a_{\ell_\nu})\\
&\geq&\ell_\nu b_{\ell_\nu}.
\end{eqnarray*}
Now let us consider the case where $k \geq \ell_{\nu+1}$. In this
case there exists $\nu'>\nu$ such that
$$k \in \{\ell_{\nu'}, \cdots \ell_{\nu'+1}-1\}.$$
Using the monotonicity of $f_1$ and $f_2$ we estimate in this case
\begin{eqnarray*}
(f_1 \cdot \sigma_* f_2)(k)&=& f_1(k)f_2(\sigma(k))\\
&=&f_1(k)f_2(\ell_{\nu'}+\ell_{\nu'+1}-k-1)\\
&\geq&f_1(\ell_{\nu}) f_2(\ell_{\nu+1}-1)\\
&=&\big(s_{\ell_{\nu}-1}(f_1)(1)\big)
\big(s_{\ell_{\nu}-1}(f_2)(\ell_{\nu+1}-\ell_\nu)\big)\\
&\geq&g_{s_{\ell_{\nu}-1}(f_1),s_{\ell_{\nu}-1}(f_2)}(\ell_{\nu+1}-\ell_\nu)\\
&\geq&\ell_\nu b_{\ell_\nu}.
\end{eqnarray*}
Hence (\ref{sigi}) and therefore Step~4 are proved.
\\ \\
\textbf{Step~5: }\emph{For every $\nu \in \N$ we have the
inequality}
\begin{equation}\label{sigi2}
\wp_{\sigma}(f_1,f_2)(\ell_\nu) \geq \ell_\nu b_{\ell_\nu}.
\end{equation}
\\
We assume by contradiction that
\begin{equation}\label{nosi}
\wp_{\sigma}(f_1,f_2)(\ell_\nu) < \ell_\nu b_{\ell_\nu}.
\end{equation}
By construction of $\wp_\sigma$ there exists $\sigma' \in
\mathfrak{S}$ such that
\begin{equation}\label{nosi2}
\wp_{\sigma}(f_1,f_2)=\sigma'_*(f_1 \cdot \sigma_* f_2) \in
\widetilde{\mathcal{F}}.
\end{equation}
Since $\wp_\sigma(f_1,f_2)$ is monotone, we deduce from (\ref{nosi})
\begin{equation}\label{nosi3}
\wp_\sigma(f_1,f_2)(k)<\ell_\nu b_{\ell_\nu}, \quad \forall\,\, k
\in \{1,\cdots, \ell_\nu\}.
\end{equation}
We define
$$A=\sigma'\big(\{1,\ldots,\ell_\nu\}\big) \subset \N.$$
By (\ref{nosi2}) and (\ref{nosi3}) we conclude
\begin{equation}\label{nosi4}
(f_1 \cdot \sigma_* f_2)(k)<\ell_\nu b_{\ell_\nu}, \quad
\forall\,\,k \in A.
\end{equation}
Since $\sigma'$ is a bijection we have
$$\#A=\ell_\nu$$
and hence it follows from (\ref{nosi4}) that there exists
$$k_0 \geq \ell_\nu$$
with the property
$$(f_1 \cdot \sigma_* f_2)(k_0)<\ell_\nu b_{\ell_\nu}.$$
But this contradicts Step~4 and hence Step~5 follows.
\\ \\
\textbf{Step~6: } \emph{The set $\mathfrak{S}_0 \cup \{\sigma\}$ is
$(f_1,f_2)$-wild.}
\\ \\
Since $\mathfrak{S}_0$ is already $(f_1,f_2)$-wild by assumption we
are left with showing that
\begin{equation}\label{wil}
[\wp_\sigma(f_1,f_2)] \neq [\wp_{\sigma'}(f_1,f_2)], \quad
\forall\,\,\sigma' \in \mathfrak{S}_0.
\end{equation}
We assume by contradiction that there exists $\sigma' \in
\mathfrak{S}_0$ such that
$$[\wp_\sigma(f_1,f_2)] = [\wp_{\sigma'}(f_1,f_2)].$$
Hence there exists $c>0$ such that
\begin{equation}\label{wi}
\wp_\sigma(f_1,f_2)(n) \leq c\big(\wp_{\sigma'}(f_1,f_2)(n)\big),
\quad \forall\,\,n \in \N.
\end{equation}
Now choose $\nu \in \N$ satisfying $\nu>c$. We estimate using Step~5
$$\wp_\sigma(f_1,f_2)(\ell_\nu)\geq \ell_\nu b_{\ell_\nu}
\geq \nu \wp_{\sigma'}(f_1,f_2)(\ell_\nu)> c
\wp_{\sigma'}(f_1,f_2)(\ell_\nu).$$ This contradicts (\ref{wi}) and
hence (\ref{wil}) has to hold. This finishes the proof of Step~6 and
hence of the Lemma. \hfill $\square$
\\ \\
\textbf{Proof of Proposition~\ref{wilder}: } We define inductively
$(f_1,f_2)$-wild subsets $\mathfrak{S}^n_0 \subset \mathfrak{S}$ of
cardinality $n \in \N$ in the following way. We set
$$\mathfrak{S}_0^1=\{\mathrm{id}\}.$$
Given $\mathfrak{S}_0^n$ there exists by Lemma~\ref{willi} $\sigma
\in \mathfrak{S} \setminus \mathfrak{S}^n_0$ such that
$\mathfrak{S}_0^n \cup \{\sigma\}$ is $(f_1,f_2)$-wild. We put
$$\mathfrak{S}^{n+1}_0=\mathfrak{S}^n_0 \cup \{\sigma\}.$$
The sets $\{\mathfrak{S}^n_0\}_{n \in \N}$ build a nested sequence
\begin{equation}\label{nest}
\mathfrak{S}^1_0 \subset \mathfrak{S}^2_0 \subset \mathfrak{S}^3_0
\subset \cdots.
\end{equation}
We define
$$\mathfrak{S}_0=\bigcup_{k=1}^\infty \mathfrak{S}^k_0 \subset
\mathfrak{S}.$$ The set $\mathfrak{S}_0$ has infinite cardinality.
We claim that it is still an $(f_1,f_2)$-wild subset of
$\mathfrak{S}$. Pick $\sigma, \sigma' \in \mathfrak{S}_0$. There
exist $j, j' \in \N$ such that
$$\sigma \in \mathfrak{S}_0^{j}, \quad \sigma'
\in \mathfrak{S}_0^{j'}.$$ We set
$$i=\max\{j,j'\}.$$
It follows from (\ref{nest}) that
$$\sigma, \sigma' \in \mathfrak{S}^i_0.$$
But since $\mathfrak{S}^i_0$ is an $(f_1,f_2)$-wild subset of
$\mathfrak{S}$ we deduce that
$$[\wp_\sigma(f_1,f_2)] \neq [\wp_{\sigma'}(f_1,f_2)].$$
This proves that $\mathfrak{S}_0$ is $(f_1,f_2)$-wild and hence we
have constructed an $(f_1,f_2)$-wild subset of infinite cardinality.
This finishes the proof of the Proposition. \hfill $\square$
\\ \\
\textbf{Proof of Theorem B: } For given $\phi_1,\phi_2 \in
\mathcal{F}$ we first choose representatives $f_1,f_2 \in
\widetilde{\mathcal{F}}$ such that
$$[f_i]=\phi_i, \quad i \in \{1,2\}.$$
By Proposition~\ref{wilder} there exists an $(f_1,f_2)$-wild subset
$\mathfrak{S}_0 \subset \mathfrak{S}$ of infinite cardinality. We
pick $\sigma \in \mathfrak{S}_0$ and introduce the triple
$$\mathcal{H}=(\ell^2,\ell^2_{f_1},\ell^2_{f_1 \cdot \sigma_*
f_2}).$$ Although $f_1 \cdot \sigma_* f_2$ is only in $\mathcal{U}$
and not necessarily in $\widetilde{\mathcal{F}}$ we define
$\ell^2_{f_1 \cdot \sigma_* f_2}$ as a subset of $\ell^2$ in the
same way as we do it in the monotone case. We further note that
$\ell^2_{f_1 \cdot \sigma_* f_2} \subset \ell^2_{f_1}$. Let
$$I \colon \ell^2_{f_1} \to \ell^2$$
be the canonical isometry as explained in the proof of
Corollary~\ref{a2}. We note that the restriction of $I$ to
$\ell^2_{f_1 \cdot \sigma_* f_2}$ gives rise to an isometry
$$I|_{\ell^2_{f_1 \cdot \sigma_* f_2}} \colon
\ell^2_{f_1 \cdot \sigma_* f_2} \to \ell^2_{\sigma_* f_2}.$$ Define
a further isometry
$$J_\sigma \colon \ell^2 \to \ell^2$$
which is given on standard basis vectors $\{\varepsilon_n\}_{n\in
\N}$ of $\ell^2$ by the formula
$$J_\sigma(\varepsilon_n)=\varepsilon_{\sigma(n)}.$$
We note that the restriction of $J_\sigma$ to $\ell^2_{\sigma_*
f_2}$ gives an isometry
$$J_\sigma|_{\ell^2_{\sigma_* f_2}} \colon \ell^2_{\sigma_* f_2}
\to \ell^2_{f_2}.$$ We deduce that the composition of $I$ and
$J_\sigma$ gives an isometry of pairs
\begin{equation}\label{supi}
J_\sigma \circ I \colon (\ell^2_{f_1},\ell^2_{f_1 \cdot \sigma_*
f_2}) \to (\ell^2,\ell^2_{f_2}).
\end{equation}
As a first consequence of the isometry (\ref{supi}) and the fact
that $(\ell^2,\ell^2_{f_2})$ is a scale Hilbert pair we conclude
that $(\ell^2_{f_1},\ell^2_{f_1 \cdot \sigma_* f_2})$ is also a
scale Hilbert pair. Since $(\ell^2,\ell^2_{f_1})$ is a further scale
Hilbert pair, we deduce that $\mathcal{H}$ is a scale Hilbert
triple. As a second consequence of (\ref{supi}) we obtain the
formula
\begin{equation}\label{in1}
\mathfrak{K}([\mathcal{H}])(1,2)=[f_2]=\phi_2.
\end{equation}
By construction of $\wp_\sigma(f_1,f_2)$ there exists $\sigma' \in
\mathfrak{S}$ such that
$$\wp_\sigma(f_1,f_2)=\sigma'_*(f_1 \cdot \sigma_* f_2).$$
Hence we obtain a scale isometry of scale Hilbert pairs
$$J_{\sigma'} \colon (\ell^2,\ell^2_{\wp_\sigma(f_1,f_2)})
\to (\ell^2,\ell^2_{f_1 \cdot \sigma_*f_2})$$ from which we deduce
\begin{equation}\label{in2}
\mathfrak{K}([\mathcal{H}])(0,2)=[\wp_\sigma(f_1,f_2)].
\end{equation}
Furthermore,
\begin{equation}\label{in3}
\mathfrak{K}([\mathcal{H}])(0,1)=[f_1]=\phi_1.
\end{equation}
Combining (\ref{in1}), (\ref{in2}), and (\ref{in3}) we obtain
$$\mathfrak{K}([\mathcal{H}])=(\phi_1,\phi_2,[\wp_\sigma(f_1,f_2)])$$
implying that
$$[\wp_\sigma(f_1,f_2)] \in \mathcal{B}(\phi_1,\phi_2).$$
Hence we get a map
$$\mathfrak{S}_0 \to \mathcal{B}(\phi_1,\phi_2), \quad
\sigma \mapsto [\wp_\sigma(f_1,f_2)]$$ which by definition of
$(f_1,f_2)$-wild is injective. Since $\#\mathfrak{S}_0=\infty$ we
deduce that
$$\#\mathcal{B}(\phi_1,\phi_2)=\infty.$$
This finishes the proof of Theorem~B. \hfill $\square$
\\ \\
\textbf{Proof of Corollary~\ref{b1}: }We only prove that the map
$\mathfrak{J}\colon \mathfrak{F} \to \mathscr{S}$ is not surjective.
The proof that $\mathfrak{J}_n \colon \mathfrak{F}_n \to
\mathscr{S}_n$ is not surjective for $n \geq 3$ is analogous, but we
prefer to avoid keeping track of the subscript $n$.
\\
By Theorem~B there exists a scale Hilbert triple
$\mathcal{H}=(\mathcal{H}_0,\mathcal{H}_1,\mathcal{H}_2)$ such that
$$\mathfrak{K}([\mathcal{H}])(0,2)
\neq \mathfrak{K}([\mathcal{H}])(0,1) \cdot
\mathfrak{K}([\mathcal{H}])(1,2).$$ Choose an arbitrary scale
Hilbert space $\mathcal{H}'=(\mathcal{H}'_0,\mathcal{H}'_1,\ldots)$.
By Corollary~\ref{a1} the Hilbert spaces $\mathcal{H}_2$ and
$\mathcal{H}'_0$ are isometric to $\ell^2$ and in particular
isometric to each other. Hence let
$$I \colon \mathcal{H}'_0 \to \mathcal{H}_2$$
be an isometry of Hilbert spaces. We now define a new scale Hilbert
space $\widetilde{\mathcal{H}}=(\widetilde{\mathcal{H}}_0,
\widetilde{\mathcal{H}}_1,\ldots)$ by setting
$$\widetilde{\mathcal{H}}_k=\left\{\begin{array}{cc}
\mathcal{H}_k & 0 \leq k \leq 2 \\
I(\mathcal{H}'_{k-2}) & k \geq 3.
\end{array}\right.$$
For the scale Hilbert space $\widetilde{\mathcal{H}}$ we still have
$$\mathfrak{K}([\widetilde{\mathcal{H}}])(0,2)
\neq \mathfrak{K}([\widetilde{\mathcal{H}}])(0,1) \cdot
\mathfrak{K}([\widetilde{\mathcal{H}}])(1,2).$$ On the other hand if
$\Phi \in \mathfrak{F}$, then we necessarily have
$$\mathfrak{K}(\mathfrak{J}(\Phi))(0,2)
= \mathfrak{K}(\mathfrak{J}(\Phi))(0,1) \cdot
\mathfrak{K}(\mathfrak{J}(\Phi))(1,2).$$ This shows that
$[\widetilde{\mathcal{H}}]$ cannot lie in the image of
$\mathfrak{J}$ and hence the Corollary is proved. \hfill $\square$

\end{document}